\documentstyle[amssymb,12pt]{amsart}

\def\set#1{\{\,#1\,\}}
\def\seq#1{\langle\,#1\,\rangle}
\def\forces{\mathrel\|\joinrel\relbar} 
\def\rest{{\restriction}}

\def\sse{\subseteq}
\def\a{\alpha}
\def\b{\beta}
\def\c{{\frak c}}  
\def\d{\delta}
\def\f{\varphi}
\def\g{\gamma}
\def\h{\eta}

\def\r{\rho}
\def\s{\sigma}
\def\w{\omega}
\def\x{\xi}

\def\A{\aleph}  

\outer\def\case #1. {\medbreak\noindent{\sl Case #1.\quad}}

\title{Ultrafilters on $\w$}
\author{James E. Baumgartner}
\address {Department of Mathematics and Computer Science,
Dartmouth College, Han\-over, NH 03755 USA}
\email {james.baumgartner@@dartmouth.edu}
\thanks{Research supported by NSF grant DMS--9113359}

\newtheorem{theorem}{Theorem}[section]
\newtheorem{corollary}[theorem]{Corollary}
\newtheorem{proposition}[theorem]{Proposition}
\newtheorem{lemma}[theorem]{Lemma}
\newtheorem{problem}{Problem}

\newcommand{\thm}[2]{\begin{theorem}\label{#1}$\!\!\!$
#2\end{theorem}}
\newcommand{\cor}[2]{\begin{corollary}\label{#1}$\!\!\!$
#2\end{corollary}}
\newcommand{\prop}[2]{\begin{proposition}\label{#1}$\!\!\!$
#2\end{proposition}}
\newcommand{\lem}[2]{\begin{lemma}\label{#1}$\!\!\!$
#2\end{lemma}}

\newcommand{\Q}{\mbox{$\Bbb Q$}}
\newcommand{\R}{\mbox{$\Bbb R$}}
\newcommand{\e}{\mbox{$\varepsilon$}}
\newcommand{\tp}{\operatorname{tp}}

     \newcommand{\charf}[1]{\mbox{\raise.48ex\hbox{$\chi$}$_{#1}$}}

\date{}

\begin{document}\maketitle
\begin{abstract}
We study the $I$-ultrafilters on $\w$, where $I$ is a a collection of subsets
of
a set $X$, usually \R\ or $\w_1$.  The $I$-ultrafilters usually contain the
$P$-points, often as a small proper subset.  We study relations between
$I$-ultrafilters for various $I$, and closure of $I$-ultrafilters under
ultrafilter sums.  We consider, but do not settle, the question whether
$I$-ultrafilters always exist.
\end{abstract}
\advance \baselineskip by 2\jot

\section{$I$-ultrafilters}

Let $I$ be a family of subsets of a set $X$ such that $I$ contains all
singletons
and is closed under subsets.  Given an ultrafilter $u$ on $\w$, we say that
$u$ is an $I$-{\em ultrafilter} if for any $F:\w\to X$ there is $A\in u$ such
that $F(A)\in I$.  Note that $X$ is recoverable from $I$ since $X=\bigcup I$. 
As values of $X$ we mostly consider ${}^\w2$, the reals \R, ${}^\w\w$, the
rationals \Q, and $\w_1$.  One may replace $\w$ in this definition by an
arbitrary set as well, but we have no use for this generalization.

If $I\sse J$ then any $I$-ultrafilter is a $J$-ultrafilter.

\prop{prop1} {If $I$ is as above then the $I$-ultrafilters are closed
downward
under the Rudin-Keisler ordering $\le_{\text{RK}}$.}

For the case where $X={}^\w2$ (or $X$ is any other version of the real
numbers), if $I$ consists of all discrete sets then we refer to the
$I$-ultrafilters as the discrete ultrafilters.  Similarly, if $I$ is the
collection of all scattered sets (or all nowhere dense sets) we refer to the
$I$-ultrafilters as the scattered ultrafilters (or the nowhere dense
ultrafilters).  If $I$ consists of all sets $Y$ with closure $\overline Y$ of
measure zero, we refer to the $I$-ultrafilters as the measure zero
ultrafilters.  Note that these definitions give the same classes of
ultrafilters when ${}^\w2$ is replaced by ${}^\w\w$ or $\Bbb R$ (or even the
interval $[0,1]$).

Observe that if $I$ is the collection of finite sets, then the
$I$-ultrafilters are exactly the principal ultrafilters.

\thm{thm2}{Suppose $u$ is an ultrafilter on $\w$.  In the following list of
conditions on $u$, each implies the
next.\begin{list}{(\roman{enumi})}{\usecounter{enumi}}\item $u$ is $P$-point
\item
$u$ is discrete\item$u$ is scattered\item $u$ is measure zero\item $u$ is
nowhere
dense\end{list}}

\begin{pf} Any discrete set is scattered, any scattered set has closure
of measure zero, and any set with closure of measure zero must be
nowhere dense, so all of this is clear except (i)$\Rightarrow$(ii), for we
have
as yet no characterization of the $P$-points as $I$-ultrafilters for some
$I$.
The following lemma will complete the proof.

\lem{lem3}{(Booth\cite{Booth})  If $I=\{Y\sse{}^\w2: Y\text{ is finite or has
order type }\w\text{ or }\w^*\}$ then the
$I$-ultrafilters are exactly the $P$-points.}
\begin{pf} Note that the order type of $Y$ is calculated with respect to
the lexicographic ordering.  Recall that an ultrafilter
$u$ is
$P$-point iff for any sequence $\seq{A_n:n<\w}$ of elements of $u$ there is
$A\in
u$ such that $A\sse^*A_n$ for all $n$.  Suppose now $u$ is
$P$-point and $F:\w\to{}^\w2$.  One easily sees by induction on
$n$ that there is $f\in{}^\w2$ such that $F^{-1}([f\rest n])\in
u$ for all $n$.  Since $u$ is $P$-point there is $A\in u$ such
that $A\sse^*F^{-1}([f\rest n])$ for all $n$  It is clear that
the unique limit point of $F(A)$ is $f$.  If $B=\set{m\in
A:F(m)<f}$ and $C=\set{m\in A:F(m)>f}$, then $F(B)\in I$,
$F(C)\in I$ and either $B\in u$, $C\in u$ or $F^{-1}\{f\}\in u$. 
Thus $u$ is an $I$-ultrafilter.

Now suppose $u$ is an $I$-ultrafilter and $\seq{A_n:n<\w}$ is a sequence of
elements of $u$.  Without loss of generality, assume $A_n\supseteq A_{n+1}$
for
all $n$ and that $\bigcap\set{A_n:n<\w}=0$.  Fix $f\in{}^\w2$ and define
$F:\w\to{}^\w2$ such that $F$ is one-to-one and if $n\in\w$ and $m$ is
maximal
with $n\in A_m$ then $F(n)\in[f\rest m]-[f\rest m+1]$.  If now $A\in u$ and
$F(A)\in I$, it is clear that $F(A)\cap[f\rest n]\ne 0$ for all $n$, so
$f$ must be a limit point of $F(A)$, and hence
$A\sse^*A_n$ for all $n$.\end{pf}\end{pf}

We are interested in the existence of ultrafilters of the types mentioned in
Theorem \ref{thm2}, but this may be an independence question. 
Shelah\cite{Ppt}
has shown it consistent that there are no $P$-point ultrafilters, and it is
conceivable that the same thing may happen with the other types of
ultrafilters
as well.  For the case of nowhere dense ultrafilters, this existence question
has already appeared as Question 31 in \cite{HartvM}, which arises from
an earlier question of van Douwen in \cite{rempt}.

With the right set-theoretic hypotheses the situation is very nice.

\thm{thm4} { Assume Martin's Axiom for $\s$-centered partial orderings.  Then
none of the implications in Theorem \ref{thm2} reverses.}

{\em Proof.} 
The first two implications are easily seen not to reverse.  Let us
apply some results from the next section.  By Theorem \ref{2x} if $u$ is
$P$-point and $n<\w$ then $u^n$ is discrete, so for $n>1$ we obtain a
discrete,
non-$P$-point ultrafilter.  Also, if $u$ is any scattered ultrafilter, then
by
Theorem \ref{2y} $u^\w$ is also scattered.  And by Theorem \ref{2z} $u^\w$ is
not discrete (as long as $u$ is nonprincipal).  Thus $u^\w$ is a scattered
nondiscrete ultrafilter.  Note that we have obtained a little more than is
claimed in Theorem \ref{thm4}, namely that if there is a discrete ultrafilter
then there is a discrete ultrafilter that is not $P$-point, and if there is a
scattered ultrafilter then there is a scattered ultrafilter that is not
discrete.  Of course it is well known that Martin's Axiom for
$\s$-centered partial orderings implies the existence of $P$-points, hence
discrete and scattered ultrafilters.

For the other two implications we make direct use of Martin's Axiom.

First let us show that there is a measure zero ultrafilter that is not
scattered.  Let us represent the rationals \Q\ as a subset of ${}^\w2$,
say as $\set{f\in {}^\w2:\exists n\ \forall m>n\ f(m)=0}$.  It will suffice
to
find an ultrafilter $u$ on \Q\ such that $u$ consists of nonscattered sets
but
for all $F:\Q\to{}^\w2$ there is $B\in u$ such that $m(\overline{F(B)})=0$. 
We
will construct $u$ inductively, using the following crucial lemma.

\lem{lem5}{Suppose $X\sse[\Q]^\w$, $|X|<2^\w$, $X$ is closed under finite
intersections, and all elements of $X$ are nonscattered.  If $F:\Q\to{}^\w2$
then there is $B\in[\Q]^\w$ such that $B\cap A$ is nonscattered for all $A\in
X$ and $m(\overline{F(B)})=0$.}

\begin{pf} Recall that if $A\sse{}^\w2$ then $A$ may be written uniquely as
$A_0\cup A_1$, where $A_0$ is scattered and $A_1$ is dense in itself.  Let us
write $d(A)=A_1$.  Of course $A$ is nonscattered iff $d(A)\ne0$, so $d(A)\ne
0$ for all $A\in X$.

Fix $F:\Q\to{}^\w2$.  We will obtain $B$ by applying Martin's Axiom to a
partial ordering $P$.

We put $p\in P$ iff $p=(x_p,S_p,X_p)$ where $x_p\in[\Q]^{<\w}$,
$S_p\in[{}^{<\w}2]^{<\w}$ and $X_p\in[X]^{<\w}$, and if
$B_p=\Q-\bigcup\set{F^{-1}[s]:s\in S_p}$ then $x_p\sse B_p$ and
$d(A\cap B_p)\ne0$ for all $A\in X$.  Put $p\le q$ iff $x_p\supseteq x_q$,
$S_p\supseteq S_q$, $X_p\supseteq X_q$, and if $A\in X_q$ then
$d(A\cap B_q)\cap x_q\sse d(A\cap B_p)\cap x_p$.

If $x_p=x_q$ and $S_p=S_q$ then $p$ and $q$ are compatible, so $P$ is
$\s$-centered.

Our final goal is to find some dense sets such that if $G\sse P$ is any
filter
meeting all of them, then $B=\bigcup\set{x_p:p\in G}$ satisfies the lemma.

For $A\in X$ let $D(A)=\set{p\in P: A\in X_p \text{ and } \exists f\in x_p\ 
f\in d(A\cap B_p}$.  Then $D(A)$ is dense, for if $q\in P$ and $f\in d(A\cap
B_q)$ then $p=(x_q\cup\{f\},S_q,X_q\cup\{A\})\le q$ and $p\in D(A)$.  For
$A\in
X$ and $n<\w$ let $E(A,n)=\set{p\in P: A\in X_p\text{ and }\forall f\in
d(A\cap B_p)\cap x_p\ \exists g\in d(A\cap B_p)\cap x_p\ f\ne g\text{ and
}|f-g|<{1\over n}}$.  Here $|f-g|$ refers to the distance between $f$ and $g$
in the standard metric on ${}^\w2$.  Then $E(A,n)$ is dense as well.  Let
$q\in
P$.  For each $f\in d(A\cap B_q)\cap x_q$ choose $g_f \ne f$ such that
$g_f\in
d(A\cap B_p)$ and $|g_f-f|<{1\over n}$, and let $x_p =x_q\cup\set{g_f:f\in
d(A\cap B_q)\cap x_q}$.  Then $p=(x_p,S_q,X_q\cup\{A\})\le q$ and $p\in
E(A,n)$.  If $G$ is a filter on $P$ meeting all the $D(A)$ and $E(A,n)$, then
$B$ intersects each $A\in X$ in a nonscattered set (since
$\bigcup\set{x_p\cap
d(A): p\in G\text{ and }A\in X_p}$ is dense in itself).  What we do not
yet know is that we may choose $B$ such that $m(\overline{F(B)})=0$.  This
requires some more dense sets.

\def\e{\varepsilon}
If $G$ is a filter then let $S=\bigcup\set{S_p:p\in G}$.  It is clear that
$\overline{F(B)}\cap\bigcup\set{[s]:s\in S}=0$ so it suffices to choose $G$
so
that $\bigcup\set{[s]:s\in S}$ has measure 1.  For this it suffices to
observe
that for each $\e>0$ the following set $D_\e$ is dense:
\[ D_\e=\set{p\in P: m(\bigcup\set{[s]:s\in S_p})>1-\e}.\]
Let $q\in P$ be given.  Fix $k$ large enough so that $S_q\sse{}^{\le k}2$. 
Let
$|x_q|=m$, let $Z=\set{(f,A):f\in x_q, A\in X_q\text{ and } f\in d(A\cap
B_q)}$, and let $|Z|=l$.  Now choose $n$ so large that $n\ge k$ and
${m+l+1\over 2^n}<\e$.  For $f\in x_q$ let $s_f=F(f)\rest n$. Note that
$\Q\cap F^{-1}[s_f]\sse B_q$.  Next observe that if $f\in \overline{d(A)}$
and
$A=A_0\cup A_1$ then $f\in\overline{d(A_0)}$ or $f\in \overline{d(A_1)}$.  It
follows that if $(f,A)\in Z$ then there is $t_{fA}\supseteq F(f)\rest k$ such
that $t_{fA}\in {}^n2$ and $f\in d(A\cap F^{-1}([t_{fA}]))$.  Finally, note
that
if $A_0\cup A_1$ meets every element of $X$ in a nonscattered set then the
same
is true of $A_0$ or $A_1$.  From this it follows that there is $r\in{}^n2$
such
that $\Q\cap F^{-1}[r]\sse B_q$.  And now, if we let $x_p=x_q$,
$S_p=S_q\cup\set{s\in {}^n2:\forall t\in S_q\ [s]\cap [t]=0\text{ and } s\ne
s_f, t_{fA}, r}$ and $X_p=X_q$ then $B_p=F^{-1}(\bigcup\set{[s_f]:f\in
x_q}\cup\bigcup\set{[t_{fA}]:(f,A)\in Z}\cup [r])\cap\Q$ so 
$m(\bigcup\set{[s]:s\in S_q})=1-{m+l+1\over 2^n}>1-\e$.  Thus $p\le q$.

This now completes the proof of the Lemma.\end{pf}

By using the Lemma repeatedly we easily construct an ultrafilter $u$ on \Q\
consisting of nonscattered sets such that for all $F:\Q\to{}^\w2$ there is
$A\in u$ with $m(\overline{F(A)})=0$.

We take a similar approach to the final part of the theorem.  It will suffice
to find an ultrafilter $u$ on \Q\ such that for all $A\in u$ $m(\overline
A)>0$
and for all $F:\Q\to {}^\w2$ there is $A\in\Q$ such that $F(A)$ is nowhere
dense. Once again the construction is based on a lemma.

\lem{lem6}{Suppose $X\sse[\Q]^\w$, $|X|<2^\w$, $X$ is closed under finite
intersections, and for all $A\in X$ $m(\overline A)>0$.  If $F:\Q\to{}^\w2$
then there is $B\in[\Q]^\w$ such that $m(\overline{B\cap A})>0$ for all $A\in
X$ and $F(B)$ is nowhere dense.}

\begin{pf}Fix $F:\Q\to{}^\w2$.  Once again we obtain $B$ by applying Martin's
Axiom to a partial ordering $P$.  Put $p\in P$ iff $p=(x_p,n_p,S_p,X_p)$,
where
$x_p\in[\Q]^{<\w}$, $x_p\ne0$, $n_p<\w$, $n_p\ge1$, $S_p\in[{}^{<\w}2]^{<\w}$
and $X_p\in[X]^{<\w}$, and if
\[ B_p=\Q-\bigcup\set{[s]: s\in {}^{n_p}2\text{ and }x_p\cap
[s]=0}-\bigcup\set{F^{-1}[s]:s\in S_p}\]
then $x_p\sse B_p$ and $m(\overline{B_p\cap A})>0$ for all $A\in X$.  Set $p<
q$ iff $x_p\supseteq x_q$, $n_p> n_q$, $S_p\supseteq S_q$, $X_p\supseteq
X_q$, if $s\in{}^{n_q}2$ and $x_q\cap[s]=0$ then $x_p\cap[s]=0$, and for all
$A\in X_q$
\[ m(\overline{A\cap B_q}-\overline{A\cap
B_p})<m(\overline{A\cap B_q})\sum_{i=n_q+1}^{n_p} 2^{-i}\]
Note that transitivity of $\le$ is clear except possibly for the last clause
above.  But if $p< q< r$ then
\[m(\overline{A\cap B_r}-\overline{A\cap B_q})<m(\overline{A\cap
B_r})\sum_{i=n_r+1}^{n_q} 2^{-i}\]
and
\[m(\overline{A\cap B_q}-\overline{A\cap B_p})<m(\overline{A\cap
B_q})\sum_{i=n_q+1}^{n_p} 2^{-i}.\]
Since $B_p\sse B_q\sse B_r$ we must have $m(\overline{A\cap B_q})\le
m(\overline{A\cap B_r})$ so 
\begin{multline*}
m(\overline{A\cap B_r}-\overline{A\cap B_p})\le 
m(\overline{A\cap
B_r}-\overline{A\cap B_q}) +  m(\overline{A\cap
B_q}-\overline{A\cap B_p})\\ < m(\overline{A\cap
B_r})\sum_{i=n_r+1}^{n_p}2^{-i}
\end{multline*}
as desired.

As before, if $x_p=x_q$, $n_p=n_q$, $S_p=S_q$ then $p$ and $q$ are
compatible. 
Hence $P$ is $\s$-centered.

Once again we seek $<\c$ dense sets in $P$ such that if $G\sse P$ is a filter
meeting all of them then $B=\bigcup\set{x_p:p\in G}$ satisfies the lemma. 
For
$A\in X$ let $D(A)=\set{p\in P:A\in X_p}$.  Then $D(A)$ is dense.

Next we assert that for $s\in{}^{<\w}2$, $D_s=\set{p\in P:\exists t\in S_p\
t\supseteq s}$ is dense.  Then if $G$ meets all the $D_s$, $F(B)$ will be
nowhere dense.  Fix $q\in P$, and let $\e$ be the minimum value among all
$m(\overline{A\cap B_q})$ for $A\in X_q$.

We will find $p\le q$ such that $p\in D_s$.  We will have $X_p=X_q$,
$S_p=S_q\cup\{t\}$ for some $t\supseteq s$, $n_p=n_q+1$, and $x_p\supseteq
x_q$
chosen so that for all $s'$ of length $n_p$ and $A\in X_p$ if $A\cap B_q\cap
[s']\ne0$ then $A\cap x_p\cap[s']\ne0$.  We must choose $t$ so that $p\le q$.

Choose $n$ so large that $2^{n-|s|}>|X_q|{2^{2n_q}\over\e}+1$.  Call $t$ {\em
bad} for $A\in X_q$ if $s\sse t$ and
\[m(\overline{A\cap B_q})-m(\overline{A\cap B_q -F^{-1}([t])})>{\e\over
2^{2n_q}}.\]
Note that there are no more than $2^{2n_q}\over\e$ sequences $t\in{}^n2$ that
are bad for $A$.  Let us call $t\in {}^n2$ {\em bad} if $s\sse t$ and for
some
$A\in X$ $m(\overline{A\cap B_q-F^{-1}([t])})=0$.  There is at most one $t$
that is bad.  It follows from the choice of $n$ that there is at least one
$t\in{}^n2$ which is not bad and not bad for any $A\in X_q$.  But now this
$t$
works.  Note that now $A\cap B_p=A\cap B_q-F^{-1}([t])$ for all $A\in X_q$.

Finally, for each $f\in\Q$ let $E_f=\set{p\in P:\text{either }f\in x_p\text{
or
}\exists s\in{}^{n_p}2\ f\in[s]\text{ and }[s]\cap x_p=0}$.  If we can show
that each $E_f$ is dense then we will be done, as the following argument
shows.  If $G$ meets all the $D(A)$, $D_s$, and $E_f$, then we know $F(B)$ is
nowhere dense.  Let $A\in X$.  Choose $q\in G\cap D(A)$, and let
$\e=m(\overline{A\cap B_q})$.  Let $T=\set{s:\exists p\in G\
s\in{}^{n_p}2\text{
and }[s]\cap x_p\ne0}$.  Since $G$ meets all the $E_f$, it is the case that
$A\cap B=A-\bigcup\set{[s]:s\in T}$.  Moreover, for any finite subset $T'\sse
T$ there is $p\in G$, $p\le q$, such that $\forall s\in T'$ $[s]\cap(A\cap
B_p)=0$.  Thus $A-\bigcup\set{[s]:s\in T'}\supseteq A\cap B_p$ and
\[m(\overline{A\cap B_q})-m(\overline{A\cap B_p})\le
\left(\sum_{i=n_q+1}^{n_p}
2^{-2i}\right)m(\overline{A\cap B_q})<{\e\over 2}.\]
Since now $\overline{A\cap B}=\bigcap_{p\in G}\overline{A\cap B_p}$ we must
have $m(\overline{A\cap B})\ge{\e\over 2}>0$ so we will be done.

Thus it suffices to prove that $E_f$ is dense.  Fix $q\in P$.  If $f\in B_q$
then put $x_p=x_q\cup\{f\}$, $n_p=n_q$, $S_p=S_q$, and $X_p=X_q$.  Then $p\le
q$ and $p\in E_f$.  So assume $f\notin B_q$.  If
$f\notin\bigcup\set{F^{-1}([s]):s\in S_q}$ then $q\in E_f$ so we may assume
$f\in F^{-1}([s])$ for some $s\in S_q$.  Find $n$ so large that $n>n_q$ and
$[f\rest n]\cap x_q=0$, and set $n_p=n$, $S_p=S_q$, $X_p=X_q$, and choose
$x_p\supseteq x_q$ so that for all $s'$ of length $n_p$ and $A\in X_p$, if
$s'\ne f\rest n$ and $A\cap B_q\cap[s']=0$ then $A\cap x_p\cap[s']=0$.  Then
$p\le q$ and $p\forces f\rest n\in T$.  This completes the proof of the
theorem.
\end{pf}

\section{Closure under sums}

If the elements of $\{u\}\cup\set{v_n:n<\w}$ are ultrafilters on $\w$ then
$\sum_u\seq{v_n:n<\w}$ is the ultrafilter on $\w\times\w$ defined by $X\in
\sum_u\seq{v_n:n<\w}$ iff $\set{n:\set{m: (n,m)\in X}\in v_n}\in u$.  We
often
identify isomorphic ultrafilters so we occasionally regard
$\sum_u\seq{v_n:n<\w}$ as an ultrafilter on $\w$.  All this also makes sense
for principal ultrafilters as well as nonprincipal ones.  We generally ignore
the principal case since it is almost always trivial, but we might note that
if
$\set{n: v_n\text{ is principal}}\in u$ then $\sum_u\seq{v_n:n<\w}$ is
isomorphic to $u$.

For $n<\w$ define $u^n$ by letting $u^0$ be principal and letting
$u^{m+1}=u\times u^m$ ($=\sum_u\seq{v_n:n<\w}$ where each $v_n=u^m$).  It is
well known and easy to see that $u^k\times u^l$ and $u^l\times u^k$ are
isomorphic for all finite $k$ and $l$.  Let $u^\w=\sum_u\seq{u^n:n<\w}$. 
This
definition may be carried out for other countable ordinals if desired.

Let $C$ and $D$ be classes of ultrafilters.  We say $C$ is {\em closed under}
$D$-{\em sums} provided that whenever $\set{v_n:n<\w}\sse C$ and $u\in D$
then
$\sum_u\seq{v_n:n<\w} \in C$.  In practice we often talk about closure under
Ramsey sums, $P$-point sums, scattered sums, etc.

If the elements of $\{u\}\cup\set{v_n:n<\w}$ are nonprincipal ultrafilters
then
$\sum_u\seq{v_n:n<\w}$ cannot be $P$-point, so closure under sums is not a
part
of the theory of $P$-points.

A family $I$ of subsets of \R\ (or ${}^\w2$) is {\em closed under discrete
unions} if whenever $\set{X_n:n<\w}\sse I$ and $\seq{G_n:n<\w}$ are disjoint
open sets, then $\bigcup_{n<\w}(X_n\cap G_n)\in I$.  If in addition $I$ is an
ideal we refer to $I$ as a {\em discrete ideal}.  Note that if $I$ is a
family
of subsets of ${}^\w2$ then in this definition we could assume that each
$G_n$
has the form $[s_n]$ for some $s_n\in{}^{<\w}2$.

{\em Remark.} The collection of scattered sets and the collection of nowhere
dense sets are discrete ideals.  The collection of all $A$ such that
$m(\overline A)=0$ is a discrete ideal.  The collection of discrete sets is
closed under discrete unions but is not an ideal.  The collection of finite
sets is an ideal but is not closed under discrete unions.  The collection of
countable sets is a discrete ideal.

A collection $I$ of subsets of ${}^\w2$ is {\em closed under scattered
unions}
iff for any $S\sse{}^{<\w}2$ if $\set{[s]:s\in S}$ is well-founded under
$\subset$ and $X_s\in I$ for $s\in S$ then $\bigcup_{s\in S}X_s\cap[s]\in I$. 
Any collection closed under scattered unions is also closed under discrete
unions.

The rationale for this definition is the fact that $X\sse{}^\w2$ is scattered
iff there is $S\sse{}^{<\w}2$ such that $\set{[s]:s\in S}$ is well-founded
under $\subset$ and for each $s\in S$ there is $f_s\in X\cap[s]$ such that 
$X=\set{f_s:s\in S}$.  The fact may be proved using an argument of Cantor and
Bendixson.  For $X\sse{}^\w2$ let $DX$ be the set of limit points of $X$ that
lie in $X$.  Define $D^\a X$ inductively by $D^0 X=X$, $D^{\a+1}X=D(D^\a X)$,
and
if $\a$ is limit, $D^\a X=\bigcap\set{D^\b X:\b<\a}$.  For some $\a<\w_1$
$D^\a
X=D^{\a+1}X$, and $X$ is scattered iff for this $\a$ $D^\a X=0$.  For $\b<\a$
$D^\b X-D^{\b+1}X$ is discrete , so there is pairwise incomparable
$S_\b\sse{}^{<\w}2$ such that for all $s\in S_\b$\quad $[s]\cap(D^\b
X-D^{\b+1}X)$ is a singleton.  If $s\in S_\b$, $t\in S_\g$ and
$[s]\subset[t]$
then clearly $\b<\g$.  Thus $S=\bigcup\set{S_\b:\b<\a}$ is as desired.

We leave to the reader the statement of a version of closure under scattered
unions appropriate for collections of subsets of \R.

\thm{discrete:scatu}{Any discrete ideal is closed under scattered unions.}

\begin{pf} Let $I$ be a discrete ideal on ${}^\w2$.  Let $S\sse {}^{<\w}2$
and
$\seq{X_s:s\in S}$ be such that $\set{[s]:s\in S}$ is well-founded under
$\subset$ and $X_s\in I$ for all $s\in S$.  We must show that $\bigcup_{s\in
S}X_s\cap[s]\in I$.  This we do by induction on the rank of $S$, which is
defined to be the smallest ordinal $\a$ for which there is $\f:S\to\a$ such
that for all $s,t\in S$ if $[s]\subset [t]$ then $\f(s)<\f(t)$.  Let
$T=\set{t\in
S:[t]\text{ is }\subset\text{-maximal}}$.  For $t\in T$ let $S_t=\set{s\in
S:[s]\subset[t]}$.  Then $S_t$ has rank $\le\f(t)<\a$ so we know by induction
that $Y_t=\bigcup_{s\in S_t}X_s\cap[s]\in I$.  Since $I$ is an ideal,
$Z_t=Y_t\cup(X_t\cap[t])\in I$ also.  But $\set{[t]:t\in T}$ is pairwise
disjoint, so since $I$ is closed under discrete unions $\bigcup\set{Z_t:t\in
T}=\bigcup_{s\in S}X_s\cap[s]\in I$. \end{pf}

\cor{scatcor}{The smallest discrete ideal containing all singletons consists
of
the scattered sets.}

\begin{pf} By the observation above about scattered sets, for any scattered
$X$
there are $S\sse{}^\w2$ and $\seq{X_s:s\in S}$ such that $\set{[s]:s\in S}$
is
wellfounded, each $X_s$ is a singleton and $X=\bigcup_{s\in S}X_s\cap[s]$.
\end{pf}

{\em Notation.} It will help to have some uniform notation for the next few
results.  Suppose $u$ is an ultrafilter on $\w$ and $F:\w\to{}^\w2$.  Then
there is unique $f\in {}^\w2$ such that for all $n$ $F^{-1}([f\rest n])\in
u$. 
We write $f=\Phi(F,u)$. This is just the limit of the sequence $F$ along $u$. 
If $\seq{v_n:n<\w}$ is a sequence of ultrafilters and
$F:\w\times\w\to{}^\w2$,
then define $F_n:\w\to{}^\w2$ by $F_n(m)=F(n,m)$ and let $f_n=\Phi(F_n,v_n)$. 
Define $F^*:\w\to{}^\w2$ by $F^*(n)=f_n$ and let $f^*=\Phi(F^*,u)$.  We will
often use this notation without specific comment.

\thm{2y}{If $I$ is a discrete ideal (containing all singletons) then the
class of $I$-ultrafilters is closed under scattered sums.}

\begin{pf} Let $\set{v_n:n<\w}$ be a set of $I$-ultrafilters, and let $u$ be
scattered.  If $\set{n:v_n\text{ is principal}}\in u$ we are done, so assume
otherwise.  Let $F:\w\times\w\to{}^\w2$.  Since each $v_n$ is an
$I$-ultrafilter, there is $A_n\in v_n$ such that $F_n(A_n)\in I$.  We may
assume $\forall m\in A_n$ $m>n$.  Since $u$ is scattered there is $A\in u$
such
that $F^*(A)$ is scattered.  Choose $S\sse{}^{<\w}2$ such that $\set{[s]:s\in
S}$      is wellfounded under $\subset$ and $\forall s\in S$ $\exists g_s\in
[s]$
$F^*(A)=\set{g_s:s\in S}$.  Fix $s\in S$, and let $C_s=\set{n:f_n=g_s}$.  Put
$C_s=D_s\cup E_s$, where $D_s=\set{n\in C_s: \set{m:F_n(m)=g_s}\in v_n}$ and
$E_s=C_s-D_s$.  For $n\in D_s$ let $B_n=\set{m\in A_n:F_n(m)=g_s}$ and for
$n\in
E_s$ let $B_n=\set{m\in A_n:F_n(m)\ne g_s, F_n(m)\in [s]\text{ and
}F_n(m)\rest
n=g_s\rest n}$.  It will suffice to show that $F(\set{(n,m): m\in B_n, n\in
A})\in I$ since $B_n\in v_n$ for all $n\in A$.

For $s\in S$ let $X_s=F(\set{(n,m):m\in B_n, n\in C_s})$.  If we can show
that
$X_s\in I$ we will be done since $X_s\sse[s]$ by the choice of the $B_n$ and
we know that $I$ is closed under scattered unions.  Let us continue to fix
$s\in S$.  Then $F(\set{(n,m):m\in B_n, n\in D_s})=\{g_s\}$ by definition of
$D_s$, so since $I$ contains all singletons we need only show
$F(\set{(n,m):n\in B_n, n\in E_s})\in I$.

For $i\ge|s|$ let $t_i$ be the unique element of ${}^{<\w}2$ such that
$|t_i|=i+1$, $t_i\rest i=g_s\rest i$ and $t_i(i)\ne g_s(i)$.  Note that if
$m\in B_n$ and $n\in E_s$ then $F(n,m)=F_n(m)\ne g_s$ but $F_n(m)\in[s]$. 
Hence $\exists i\ge|s|$ $F_n(m)\in[t_i]$.  Note also that if $F_n(m)\in[t_i]$
then $i\ge n$ since $F_n(m)\rest n=g_s\rest n$.  Thus if $Y_i=F(\set{(n,m):
m\in
B_n, n\in E_s})\cap[t_i]$ then $Y_i\sse\set{F(B_n):n\le i}$.  Hence $Y_i\in
I$
since $F(B_n)\in I$ for all $n\le i$.  Also, since the $[t_i]$ are pairwise
disjoint and $I$ is closed under discrete unions, $\bigcup\set{Y_i:i\ge
|s|}\in
I$.  But $\bigcup\set{Y_i:i\ge |s|}=F(\set{(n,m): m\in B_n, n\in E_s}$, so
the
proof is complete. \end{pf}

{\em Remark.} Since one possibility for $I$ is the scattered sets, it is
clear
that we cannot increase the condition on $u$ independent of $I$.  It is
conceivable that in Theorem \ref{2y} the class of $I$-ultrafilters is
even closed under $I$-ultrafilter sums but we do not see how to prove this. 
However, see the next two results.  We will often apply Theorem \ref{2y}
to the case of $P$-point sums.  For example, if we begin with the $P$-points
and close under $P$-point sums, every ultrafilter so obtained must be
scattered.

\thm{nwdsum}{The collection of nowhere dense ultrafilters is closed under
nowhere dense sums.}

\begin {pf}Suppose $\{u\}\cup\set{v_n:n<\w}$ are nowhere dense ultrafilters. 
We will show that $\sum_u\seq{v_n:n<\w}$ is nowhere dense.  Suppose
$F:\w\times\w\to{}^\w2$.  Let $A_n\in v_n$ be such that $F_n(A_n)$ is nowhere
dense, and let $A\in u$ be such that $F^*(A)$ is nowhere dense.

Let $S\sse{}^{<\w}2$ be such that $\bigcup\set{[s]:s\in S}$ is dense and
disjoint from $F^*(A)$.  Let $\seq{s_n:n<\w}$ enumerate ${}^{<\w}2$.  For
each
$n$ we will find $t_n\supseteq s_n$ and, if $n\in A$, $B_n\sse A_n$ such that
$B_n\in v_n$ and $F_n(B_n)\cap[t_i]=0$ and $F_i(B_i)\cap[t_n]=0$ for all
$i\le
n$.  Proceed as follows.  First find $t\supseteq s_n$ such that $[t]\sse
\bigcup\set{[s]:s\in S}$.  Now $\bigcup\set{F_i(B_i):i<n}$ is nowhere dense
so
there is $t_n\supseteq t$ such that  $[t_n]\cap F_i(B_i)=0$ for all $i<n$.

If $n\in A$ then $f_n\notin\bigcup\set{[s]:s\in S}$ so
$f_n\notin\bigcup\set{[t_i]:i<n}$.  Let
$B_n=A_n-\bigcup\set{F_n^{-1}[t_i]:i\le n}$.  Then $B_n\in v_n$.

Thus $\bigcup\set{[t_n]:n<\w}$ is dense open and disjoint from
$\bigcup\set{F_n(B_n):n\in A}=F(\set{(n,m): m\in B_n, n\in A})$, and the
proof
is complete. \end{pf}

\thm{m0sum}{The class of measure zero ultrafilters is closed under measure
zero
ultrafilter sums.}

\begin{pf}Suppose $u$ and $v_n$, $n<\w$, are measure zero ultrafilters.  We
must show $\sum_u\seq{v_n:n<\w}$ is a measure zero ultrafilter.  Fix
$F:\w\times\w\to{}^\w2$.  Let $A_n\in v_n$ be such that
$m(\overline{F_n(A(n)})=0$, and get $A\in u$ with $m(\overline{F^*(A)})=0$. 
For $n\in A$ let $B_n=A_n\cap F_n^{-1}([f_n\rest n])$.  Then $B_n\in v_n$. 
If 
$B=\set{(n,m): n\in A, m\in B_n}$ then $B\in\sum_u\seq{v_n:n<\w}$ and
$\overline{F(B)}=\bigcup_{n\in A}\overline{F_n(B_n)}\cup\overline{F^*(A)}$,
so
$m(\overline{F(B)})=0$ as desired. \end{pf}

If $I$ is closed under discrete unions but is not an ideal, then the
situation
is not so nice.  For example, the collection of discrete ultrafilters is not
closed under $P$-point (or even Ramsey) sums, as we shall see.

\thm{2z}{For any nonprincipal ultrafilter $u$ on $\w$, $u^\w$ is not
discrete.}

\begin{pf}  Let us define a linear ordering on ${}^{<\w}\w$ by setting
$\s\le\tau$ iff either $\tau\sse\s$ or else if $i$ is minimal such that
$\s(i)\ne\tau(i)$ then $\s(i)<\tau(i)$.  Note that in the order topology $\s$
is the limit point of $\set{\s^\frown m:m\in A}$ for any infinite $A\sse\w$. 
For each $n\ge 1$ define an ultrafilter $u_n$ on ${}^n\w$ as follows.  If
$n=1$, put $X\in u_1$ iff $\set{m:\langle m\rangle\in X}\in u$.  If $n>1$ put
$X\in u_n$ iff $\set{\s:\set{m:\s^\frown m\in X}\in u}\in u_{n-1}$.  Clearly
$u_n$ is isomorphic to $u_{n-1}\times u$, so by induction $u_n$ is isomorphic
to $u^n$.  Finally, define $u_\w$ on ${}^{<\w}\w$ by $X\in u^\w$ iff
$\set{n:X\cap{}^n\w\in u_n}\in u$.  Then $u_\w$ is isomorphic to $u^\w$.

It is easy to see that if $X\in u_n$ then if $\overline X$ is the closure of
$X$ in the order topology, $\overline X\cap{}^{n-1}\w\in u_{n-1}$, and thus
by
induction $\overline X\cap{}^i\w\in u_i$ for all $i$, $1\le i<n$.  Let $X\in
u_\w$.  Choose $m,n$ such that ${}^m\w\cap X\in u_m$, ${}^n\w\cap X\in u_n$
and
$m<n$.  Then $\overline{{}^n\w\cap X}\cap{}^m\w\in u_m$ so ${}^m\w\cap X$
contains a limit point of ${}^n\w\cap X$ and $X$ is not discrete.  Of course
we
have not yet fixed a map of ${}^{<\w}\w$ into ${}^\w2$, but if
$f:{}^{<\w}\w\to{}^\w2$ is any homeomorphic embedding (relative to the order
topology) it is clear that for all $X\in u_\w$, $f(X)$ is not discrete, and
hence $u_\w$ is not a discrete ultrafilter. \end{pf}

A set $X\sse{}^\w2$ has level $\le k$ if $d^kX=0$.  If
$I=\set{X\sse{}^\w2:X\text{ has level }\le k}$ then an $I$-ultrafilter is a
level $\le k$ ultrafilter.

\thm{levelk}{If $u$ is $P$-point, $k<\w$, and for all $n<\w$ $v_n$ has level
$\le k$, then $\sum_u\seq{v_n:n<\w}$ has level $\le k+1$.}

\begin{pf} Let $F:\w\times\w\to{}^\w2$ as usual, and fix $A_n\in v_n$ with
$F_n(A_n)$ with level $\le k$.  Let $A\in u$ with $F^*(A)$ a singleton or of
order type $\w$ or $\w^*$.

{\em Case 1.} $F^*(A)=\{g\}$.  For each $n$ choose $B_n\sse A_n$ so that
$B_n\in v_n$ and $\forall m\in B_n$ $F_n(m)\in[g\rest n]$.  It is not
difficult
to see that $F\set{(n,m):m\in B_n, n\in A}$ has level $\le k+1$.

{\em Case 2.} $F^*(A)$ has order type $\w$.  (The case $\w^*$ is similar, and
will be omitted.)  Since $u$ is $P$-point, we may assume that $\forall g\in
F^*(A)$ $\set{n\in A:f_n=g}$ is finite.  If now we choose $B_n\sse A_n$ so
that
$B_n\in v_n$ and $\forall m\in B_n$ $F_n(m)\in[f_n\rest n]$ then again the
reader may check that $F\set{(n,m):m\in B_n, n\in A}$ has level $\le
k+1$. \end{pf}

\prop{levelkdiscrete}{For $k<\w$ every level $\le k$ ultrafilter is
discrete.}

\begin{pf} Suppose $F:\w\to{}^\w2$, $u$ is level $\le k$ and $A\in u$ is such
that $F(A)$ has level $\le k$.  For $i<k$ the sets
$F(A)\cap(d^iF(A)-d^{i+1}F(A))$ are disjoint and discrete, and one of them
must
have inverse image lying in $u$. \end{pf}

\cor{2x}{If $u$ is $P$-point then $u^n$ is discrete for all $n<\w$.}

\begin{pf} By Theorem \ref{levelk} each $u^n$ has level $\le n$, and by
Proposition \ref{levelkdiscrete} each is discrete. \end{pf}

\section{Forcing results}

\def\lbb{\lbrack\!\lbrack}
\def\rbb{\rbrack\!\rbrack}
In \cite{Ppt} Shelah proved it consistent that there do not exist $P$-point
ultrafilters.  Naturally, one would like to ask whether the same is true for
nowhere dense ultrafilters, and for the other classes of ultrafilters
considered so far.  We begin this section with an observation of Shelah,
presented with his permission, that shows the nonexistence of nowhere dense
ultrafilters does not follow from the nonexistence of $P$-points.

\thm{NWDnoPpt}{(Shelah) It is consistent that there are no $P$-points but
there
is a nowhere dense ultrafilter.}

{\em Proof.} Since the consistency of the nonexistence of $P$-points may be
obtained by iterating ${}^\w\w$-bounding forcing extensions with countable
support in such a way that arbitrary ${}^\w\w$-bounding extensions (like
random-real forcing) occur cofinally, it will suffice to prove the following.

\lem{NWDext}{Let $F$ be a nonprincipal filter of subsets of \Q\ (so every
element of $F$ is infinite).  If $r$ is random over $V$ then in $V[r]$ there
is
$A\sse\Q$ such that $\forall B\in F$ $A\cap B$ is infinite and $A$ is nowhere
dense.}

{\em Proof.} Let us assume that $\Q\sse{}^\w2$.  We may assume $F$ is an
ultrafilter if we wish.  Let $f\in{}^\w2$ be such that $\forall n$
$\Q\cap[f\rest n]\in F$.

Let $\seq{a_n:n<\w}$ be a decomposition of $\w$ into disjoint sets such that
$|a_n|=n+1$.  Given $g\in{}^\w2$ let $\s_n(g)\in{}^{n+1}2$ be defined as
follows.  Let $a_n=\{k_0,\ldots,k_n\}$ and let $\s_n(g)(i)=g(k_i)$ for all
$i\le n$.  If $r$ is random over $V$ then let $A=\Q-\bigcup\set{[\s_n(r)]:
n\in\w\text{ and }\s_n(r)\not\sse f}$.  We must show that $A$ works.

\lem{NWDl1}{$\forall s\in{}^{<\w}2$ $\exists n$ $\s_n(r)\supseteq s$ and
$\s_n(r)\not\sse f$.}

\begin{pf} Given $s$, we may extend it if necessary so that $s\not\sse f$. 
But
now
\[\bigwedge_{n\ge|s|}\lbb s\not\sse\s_n(r)\rbb=0\]
since the sets are all independent of measure $1-2^{-|s|}$. \end{pf}

Hence $A$ is nowhere dense.

\lem{NWDl2}{$\forall B\in F$ $A\cap B$ is infinite.}

\begin{pf} It will suffice to show that if $b$ is a condition and
$g_0,\ldots,g_k\in {}^\w2$ then $\exists g\in B$ $\exists c\le b$ $g\ne g_i$
for $i\le k$ and $c\forces g\in A$.  Let $b$ and $g_0,\ldots,g_k$ be fixed. 
Fix $n$ so that $2^{-(n-1)}<m(b)$.  Choose $g\in B\cap[f\rest n]$ such that
$g\ne
g_i$ for $i\le k$.  For $m\ge n$, note that $\lbb g\rest
m+1\ne\s_m(r)\rbb$ has measure $1-2^{m+1}$.  Thus
$c=\lbb\forall m\ge n\ g\rest m+1\ne\s_m(r)\rbb \wedge
b$ has positive measure.  But clearly $c\forces g\in A$. \end{pf}

The same trick does not work for scattered ultrafilters, so it is conceivable
that the existence of a $P$-point ultrafilter follows from the existence of a
scattered ultrafilter.

\thm{mythm}{(CH) There is an ultrafilter $u$ on \Q\ such that if $B$ is any
measure algebra then $V^B\models$ Every scattered subset of \Q\ is disjoint
from some element of $u$.}

{\em Proof.} The key ingredient of the proof is the following:

\lem{mylemma}{Let $A\sse\Q$ be dense-in-itself and suppose $\forces_B\dot X$
is
scattered.  Then $\forall\e>0$ $\exists a\in A$ $m(\lbb a\in\dot
X\rbb)<\e$.}

\begin{pf} Working in $V[G]$, let $X=\dot X_G$.  Since $X$ is scattered,
there
is a sequence $\seq{s_q:q\in X}\sse{}^\w2$ such that $q\in[s_q]$ for $q\in
X$,the $s_q$ are distinct and $\set{[s_q]:q\in X}$ is wellfounded under
$\subset$.  Since $V[G]$ is ${}^\w\w$-bounding over $V$ there is $g\in
{}^A\w\cap V$ such that $\forall q\in X\cap A$ $|s_q|\le g(q)$.  Thus we may
assume that for $q\in X\cap A$, $s_q=q\rest g(q)$.  (Note that any
counterexample to the wellfoundedness of the $[s_q]$ immediately yields a
counterexample to the wellfoundedness of the original $[s_q]$.)  We may also
assume that $g$ is one-to-one.

Of course, for $q\in A$ $q\in[q\rest g(q)]$, and since $A$ is not scattered
$\seq{[q\rest g(q)]:q\in A}$ is not wellfounded, so there is $\seq{q_n:n<\w}$
such that $q_{n+1}\rest g(q_{n+1})\supset q_n\rest g(q_n)$ for all $n$.  Let
$s_n=q_n\rest g(q_n)$.  Let $\dot S$ be a name for the set $S=\set{s_q:q\in
A\cap X}$ in $V[G]$.  For each $n$ let $b_n=\bigvee_{m\ge n}\lbb
s_n\in \dot S\rbb$.  Then clearly $\seq{b_n:n<\w}$ is a decreasing
sequence of elements of $B$.

There are two cases.  First suppose $\lim_{n\to\infty}m(b_n)=0$.  Then find
$n$
such that $m(b_n)<\e$.  Then by the property of $g$ (which we assume is
forced
by some $b\in B$, below which we work)
\[b\le\lbb\text{if }q_n\in\dot X\text{ then }s_n=q_n\rest
g(q_n)\in\dot S\rbb\]
so $b\wedge\lbb q_n\in \dot X\rbb\le b_n$ and $m(b\wedge\lbb q_n\in\dot
X\rbb)<\e$.

Now suppose $\lim_{n\to\infty}m(b_n)=\d>0$.  But now if $a=\bigwedge_n b_n$
then $m(a)=\d$ and $a\forces \set{n:s_n\in\dot S}$ is infinite, so $a\forces
\dot S$ is not wellfounded, which is impossible.

We may choose $b$ and $g$ so that $m(b)>1-{\e\over 2}$.  If we choose $n$
large
enough in Case 1 so that $m(b\wedge\lbb q_n\in\dot X\rbb)<{\e\over 2}$, then
$m(\lbb q_n\in\dot X\rbb)<\e$ and we are done. \end{pf}

Of course we could also work below any given element of $B$.

For the Theorem, note that it will suffice to prove this for $B=$ the algebra
for adding a single random real.  Let $\set{(b_\a,\dot Y_\s):\a<\w_1}$
enumerate
all $(b,\dot Y)$ such that $b\in B$ and $b\forces\dot Y$ is scattered, $\dot
Y\sse\Q$.  Now construct $\seq{A_\a:\a<\w_1}$ generating an ultrafilter on
\Q\
such that every finite intersection of the $A_\a$'s is nonscattered and
$\forall\a$ $\exists b\le b_\a$ $b\forces A_\a\cap\dot Y_\a=0$.  The latter
may
easily be accomplished by repeated use of the Lemma.

\section{Ordinal ultrafilters}

In this section we continue our study of $I$-ultrafilters, but turn our
attention to families $I$ of subsets of $\w_1$ rather than \R.  We
investigate
some classes of ultrafilters largely incomparable with those already studied.

For $\a\le\w_1$ let $I_\a=\set{X\sse\w_1:X\text{ has order type }\le\a\ (\tp
X\le\a)}$ and $J_\a=\set{X\sse\w_1: \tp X<\a}$.  Then $I_\a$ is closed under
subsets and, if $\a$ is an indecomposable ordinal, then $J_\a$ is an ideal. 
(Recall that the indecomposable ordinals are exactly the ordinal powers of
$\w$.)  For convenience we refer to $I_\a$-ultrafilters as $\a$-ultrafilters. 
A
{\em proper} $\a$-ultrafilter is an $\a$-ultrafilter that is not a
$\b$-ultrafilter for any $\b<\a$.  If $u$ is a proper $\a$-ultrafilter then
$\a$
must be indecomposable.  Note that an $\w^\a$-ultrafilter is a
$J_\b$-ultrafilter where $\b=\w^{\a+1}$.  In general we do not know whether,
if
$\a$ is limit, there is a $J_{\w^\a}$-ultrafilter that is not a
$\b$-ultrafilter for some $\b<\w^\a$, even if CH or MA is assumed.

Note that every ultrafilter on $\w$ is an $\w_1$-ultrafilter.  By Proposition
\ref{prop1} the
$\a$-ultrafilters are closed downward under the Rudin-Keisler ordering
$\le_{\text{RK}}$.

Note also that the 1-ultrafilters are exactly the principal ultrafilters. 
The
next indecomposable ordinal after $\w^0=1$ is $\w^1=\w$.

\thm{wufs}{The $\w$-ultrafilters are the $P$-point ultrafilters.}

\begin{pf} It follows from Booth's Lemma \ref{lem3} that every $P$-point is
an
$\w$-ultrafilter, for any map from $\w$ to $\w_1$ has countable range, hence
may be regarded as a map from $\w$ into ${}^\w2$.  Suppose
$u$ is an
$\w$-ultrafilter.  Let
$A_n\in u$ for
$n<\w$.  We seek $A\in u$ such that $A-A_n$ is finite for all $n$.  Without
loss of generality we may assume that $A_{n+1}\sse A_n$ and $A_n-A_{n+1}$ is
infinite for all $n$, and that $\bigcap\set{A_n:n<\w}=0$.  Choose
$f:\w\to\w^2$
one-to-one such that all the elements of $f(A_n-A_{n+1})$ precede all the
elements of
$f(A_{n+1})$ for all $n$.  If now $A\in u$ and $\tp f(A)\le\w$ then clearly
$A-A_n$ is finite for all $n$, so $u$ is a $P$-point.  \end{pf}

\thm{properas}{If there exists a $P$-point, then there are proper
$\w^\a$-ultrafilters for all $\a<\w_1$.}

\begin{pf} This is clear for $\a=0,1$, so assume $\a>1$.  We proceed by
induction.  If $\a=\b+1$ let $\a_n=\w^\b$ for all $n$ and if $\a$ is limit
choose $\a_n$, $n<\w$, indecomposable and cofinal in $\w^\a$.  For each $n$,
let
$v_n$ be a proper $\a_n$-ultrafilter and let $u$ be $P$-point.  Then we claim
$v=\sum_u\seq{v_n:n<\w}$ is a proper $\w^\a$-ultrafilter.  For each $n$
let $F_n:\w\to\w_1$ be such that for all $A\in v_n$ $\tp(F(A))\ge\a_n$, and
assume that whenever $m<n$ then the range of $F_m$ lies below the range of
$F_n$. Define
$F:\w\times\w\to\w_1$ by
$F(n,m)=F_n(m)$.  It is clear from the construction that if $A\in v$ then
$\tp F(A)\ge\w^\a$.

Now we must check that for an arbitrary
$F:\w\times\w\to\w_1$ $\exists A\in v$ 
$\tp F(A)\le\w^\a$.

For each $n$ find $A_n\in v_n$ such that $\tp F_n(A_n)\le\a_n$
(using the notation from \S 2) and let $\g_n=\sup F_n(A_n)$.  Since $u$ is an
$\w$-ultrafilter there is $A\in u$ such that $\tp\set{\g_n:n\in A}=\w$
or 1.

First suppose $\g_n=\g$ for all $n\in A$.  Then let $\seq{\d_n:n<\w}$ be
cofinal
in
$\g$ and let $B_n=A_n\cap F_n^{-1}(\g-\d_n)$.  Then $B_n\in v_n$ and
$B=\set{(n,m):m\in B_n, n\in A}$ has the property that $B\in v$ and
$\tp F(B)\le\w^\a$.

Now suppose $\tp\set{\g_n:n<\w}=\w$.  If for each $n$ we let
$C_n=\set{m:\g_m\ge\g_n}$ then $C_n\in u$ so we may assume $A-C_n$ is finite
for all $n$.  And now if $B=\set{(n,m):m\in A_n, n\in A}$ then $B\in v$
and $\tp F(B)\le\w^\a$.  \end{pf}

The construction in this theorem is rather limited, since by Theorem \ref{2y}
every ultrafilter we construct is scattered.  There are, however, other ways
to
proceed.

\thm{CHthm}{Assume Martin's Axiom for $\s$-centered partial orderings.  Then
there is an
$\w^2$-ultrafilter that is not nowhere dense.}

{\em Proof.}  Set $\Q=\bigcup\set{A_n:n<\w}$ where the $A_n$ are disjoint and
dense.  Fix $f\in{}^\w2$.  Call $A\sse\Q$ {\em good} if $\set{n:\exists
s\supseteq f\rest n\ X\cap A_n\text{ is dense in }[s]}$ is infinite.  We
construct an $\w^2$-ultrafilter on \Q\ out of good sets.

Suppose $X$ is a collection of good sets closed under finite intersection,
and
$|X|<\c$.  Let us say $A\sse\Q$ is {\em large} for $X$ if $A\cap B$ is good
for
all $B\in X$.

\lem{CHlem1}{If $A\cup B$ is large for $X$ then either $A$ or $B$ is large
for
$X$.}

\lem{CHlem2}{If $g:\Q\to\w_1$ then there is $A$ large for $X$ such that
$\tp g(A)\le\w^2$.}

\begin{pf} Let us say that $(k,s,\a)$ is {\em appropriate for} $A\in X$ iff
$k<\w$, $s\supseteq f\rest k$ and $\forall \b<\a$ $(A\cap A_k)\cap
g^{-1}(\a-\b)$ is dense in $[s]$.  If $A$ is good there are arbitrarily large
$k<\w$ for which $\exists s\supseteq f\rest k$ $A\cap A_k$ is dense in $[s]$. 
If we fix such a $k$ and let $\a$ be minimal such that $\exists s\supseteq
f\rest k$ $A\cap A_k$ is dense in $[s]$ and $\a=\sup g(A\cap A_k\cap[s])$,
then
$(k,s,\a)$ is appropriate for $A$.  Thus the set of $(k,s,\a)$ appropriate
for
$A$ is always infinite and contains triples with arbitrarily large $k$. 
Hence
by Martin's Axiom there is $h:\w\to\w\times{}^{<\w}2\times\w_1$ such that if
$h(n)=(k_n,s_n,\a_n)$ then $k_n<k_{n+1}$ and for all $A\in X$ $\exists n$
$\forall m>n$ $h(m)$ is appropriate for $A$.  (Note that if $(k,s,\a)$ is
appropriate for $A$ then we may take $\a\le\sup g(\Q)$ so only countably many
$(k,s,\a)$ need be considered).  And now we may assume that either $\forall
n$
$\a_n<\a_{n+1}$ or $\exists \a$ $\forall n$ $\a_n=\a$.

{\em Case 1.} $\a_n<\a_{n+1}$ for all $n$.  For each $n$ let
$\seq{t^n_i:i<\w}$
enumerate $\set{t\in{}^{<\w}2: t\supseteq s_n}$.  Now let $P$ consist of all
functions $p$ such that $p$ maps a finite subset of $\w\times\w$ one-one into
\Q\ and if $p(i,j)$ is defined then $p(i,j)\in A_i\cap[t^i_j]$, $\a_{i-1}\le
g(p(i,j)) <\a_i$, and $g(p(i,j))\le g(p(i,j+1))$.  Put $p\le q$ iff
$p\supseteq
q$.  Then $P$ is $\s$-centered so we may apply Martin's Axiom to it.  Clearly
$D_n=\set{p:n\times n\sse\operatorname{domain}(p)}$ is dense, and if $A\in X$
and $n<\w$ then $D(A,n)=\set{p:\forall i<n\text{ if }(k_i,s_i,\a_i)\text{ is
appropriate for }A\text{ then }\set{p(i,j):j\text{ is such that }p(i,j)\text{
is defined}}\cap A\text{ has cardinality }\ge n}$ is dense as well.  If $G$
is
generic for all the $D_n$ and $D(A,n)$ and
$B=\bigcup\set{\operatorname{range}(p): p\in G}$ then $B\cap A$ is good for
all
$A\in X$ and $\tp g(B)\le\w^2$.

{\em Case 2.} $\a_n=\a$ for all $n$.  If $\a=\b+1$ then $g^{-1}\{\b\}\cap A$
is
good for all $A\in F$, so we are done.  Assume $\a$ is limit.  Define the
$t^n_i$ as in Case 1 and define $P$ as we did there as well, except that we
drop the condition $\a_{i-1}\le g(p(i,j))$ and add the requirement that
$g(p(i-1,j))\le g(p(i,j))$ for all $i\ge 1$ and all $j$.  Then we are done as
before. \end{pf}

Using the lemmas it is straightforward to build an $\w^2$-ultrafilter out of
good sets.  Such an ultrafilter cannot be nowhere dense.

The $P$-point ultrafilters, however, cannot be entirely eliminated.

\thm{needppt}{Let $\a<\w_1$ and assume $u$ is a proper
$\w^{\a+1}$-ultrafilter.  Then there is a $P$-point ultrafilter $v$ such that
$v\le_{\text{RK}}u$.}

\begin{pf} We may assume $u$ is an ultrafilter on $\w^{\a+1}$ such that
$\forall A\in u$ $\tp A=\w^{\a+1}$.  Now define
$f:\w^{\a+1}\to\w$ by $f(\x)=n$ iff $\w^\a\cdot n\le\x<\w^\a\cdot(n+1)$.  Let
$v=f(u)$.  If $v$ is not an $\w$-ultrafilter then there is $g:\w\to\b$ such
that
$\forall B\in v$ $\tp g(B)>\w$.  Define
$h:\w^{\a+1}\to\w^\a\cdot\b$ by $h(\x)=
\w^\a\cdot g(f(\x)))+(\x-\w^\a\cdot f(\x))$  (Here by $\x-\w^\a\cdot f(\x)$
we
mean that $\h$ such that $\w^\a\cdot f(\x)+\h=\x$.)  Let $A\in u$.  Clearly
we
must have $B=\set{n:\tp(A\cap f^{-1}\{n\})=\w^\a}\in v$.  But since $\tp
g(B)>\w$ we must have $\tp h(a)>\w^\a\cdot\w=\w^{\a+1}$, a contradiction.
\end{pf}

\cor{pptcor}{If there are no $P$-points, then there are no proper
$\w^{\a+1}$-ultrafilters for $\a<\w_1$.}

{\em Remark.} Several questions come to mind.  Is it consistent that there
are
no $P$-points but $\w^\w$-ultrafilters exist?  Assuming CH or MA, does it
follow that there exists an $\w^\w$-ultrafilter $u$ such that $\forall
v\le_{\text{RK}}u$ if $v$ is nonprincipal then $v$ is a
proper $\w^\w$-ultrafilter?  If $u$ is a proper $\w^{\a+\w}$-ultrafilter,
does it
follow that there is a proper $\w^\w$-ultrafilter $v\le_{\text{RK}}u$?

It is clear from the proof of Theorem \ref{properas} that if we begin with
the
$P$-points and close under $P$-point sums, we will obtain proper
$\a$-ultrafilters for arbitrarily large $\a<\w_1$.  Thus the following result
is about the best we can hope for.

\thm{asums}{Let $C=\set{u:\exists\a<\w_1\ u\text{ is a $J_\a$-ultrafilter}}$. 
Then $C$ is closed under $C$-sums.}

\begin{pf}Assume $\set{v_n:n<\w}\cup\{u\}$ are $J_\a$-ultrafilters.  It will
suffice to show that $\sum_u\seq{v_n:n<\w}$ is a $J_\r$-ultrafilter, where
$\r=\a\cdot\w^\a$.  We may assume that $\a$ is indecomposable.  Let
$F:\w\times\w\to\w_1$ and let $F_n(m)=F(n,m)$ as before.  Fix $A_n\in v_n$
with
$\tp F_n(A_n)=\a_n$ minimal.  Then $\a_n$ is indecomposable, so is either a
limit ordinal or is 1.  Define $F^*(n)=\sup F_n(A_n)$ and fix $A\in u$ with
$\tp F^*(A) <\a$.

If $\set{n:\a_n=1}\in u$ then it is easy to find a set in
$\sum_u\seq{v_n:n<\w}$ whose $F$-image has order type $<\a$ (just consider
$\bigcup\set{A_n:\a_n=1}$), so we may assume that $\a_n$ is a limit ordinal
for
all $n\in A$.  Let $G\supseteq F^*(A)$ be a countable subset of $\w_1$ such
that every element of $F^*(A)$ is a limit point of $G$, and let $\pi:\w\to G$
be a bijection.  For each $n$, let $B_n=\set{m\in A_n:\forall i\le n\text{ if
}\pi(i)<F^*(n)\text{ then }\pi(i)<F(n,m)}$.  Then $B_n\in v_n$.  We claim
$F(\bigcup\set{B_n:n\in A})$ has type $<\a\cdot\w^\a$.

For $\g\in F^*(A)$ let $B_\g=\bigcup\set{B_n:F^*(n)=\g}$.  Note that if
$\g'<\g$ and $\g'\in G$, say $\pi(i)=\g'$, then $F(B_\g)\cap\g'\sse
F(\bigcup\set{B_j:j\le i})$ which has order type $<\a$, so $F(B_\g)$ has
order
type $\le\a$.

Now by induction on $\b<\a$ we prove that $F(\bigcup\set{B_{\g'}:\g'<\g,
\g'\in
F^*(A)})$ has type $\le \a\cdot\w^\b$, where $\g$ is the $\b$th element of
the
closure of
$F^*(A)$.  This is clear if $\b=0$ or if $\b$ is a successor ordinal (since
$\a\cdot\w^\b$ is always indecomposable).  So suppose $\g$ is a limit point
of
$F^*(A)$.  Note that if $\g'<\g$, $\g'\in F^*(A)$, and
$n(\g')=\min\set{n:F^*(n)=\g'}$ then the $n(\g')$ are cofinal in $\w$, and
that
$F^*(B_\g')\cap\pi(i)=0$ for all $i<n(\g')$ with $\pi(i)<\g'$.  Hence by
inductive hypothesis $F(\bigcup\set{B_{\g'}:\g'<\g, \g'\in F^*(A)})$ has type
$\le\a\cdot\w^\b$, and we are done. \end{pf}

Let us finish this section by considering ultrafilters that are not
$J_\a$-ultrafilters for any $\a<\w_1$.

\thm{zfcthm}{(ZFC) There is an ultrafilter that is not a $J_\a$-ultrafilter
for
any $\a<\w_1$. ($\a>1$, $\a$ indecomposable.)}

{\em Proof.} The referee has commented that there is an easy proof, which
goes
as follows.  For any $\a<\w_1$ there is a $u_\a$ that is not a
$J_\a$-ultrafilter; simply choose indecomposable $\b>\a$ and let $u_\a$ be an
ultrafilter on $\b$ disjoint from $J_\b$.  It is well known (see Comfort and
Negrepontis \cite{Comfort}, Corollary 10.13(a)) that for any collection of
$\le
2^{\A_0}$ ultrafilters on
$\w$ there is an ultrafilter $u$ lying above all of them in the Rudin-Keisler
ordering.  Of course $u$ cannot be a $J_\b$-ultrafilter for any
indecomposable
$\b<\w_1$ by Proposition \ref{prop1}.

We will now, however, give another simple proof of this result that is
relevant
to the proof of Theorem \ref{incomp}.

Let $\a<\w_1$.  An $\w$-{\em block of} $\a$ is any set of the form
$\set{\x:\b\le\x<\b+\w}$ where $\b+\w\le\a$.

\lem{zfclem1}{For each indecomposable $\a<\w_1$ there is $f_\a:\w\to\a$, a
surjection such that whenever $\a_0<\a_1<\ldots<\a_n=\a$ are indecomposable
and
$A_i$ is an $\w$-block of $\a_i$ then $f_{\a_0}^{-1}A_0\cap\ldots\cap
f_{\a_n}^{-1}A_n$ is infinite.}

\begin{pf} Simply take independent maps, or construct the $f_\a$ by
induction on $\a$.
\end{pf}

\lem{zfclem2}{If $\a$ is indecomposable, $X\sse\a$ and $\tp X<\a$ then there
is
an $\w$-block $A$ of $\a$ (in fact many disjoint such blocks) such that
$A\cap
X=0$.}

Now let $I$ be the ideal generated by all sets of the form $f_\a^{-1}X$,
where
$\a$ is indecomposable, $X\sse \a$ and $\tp X<\a$.

\lem{zfclem3}{$I$ is proper.}

\begin{pf} Suppose $X_i\sse\a_i$, $\tp X_i<\a_i$ for $i< n$.  Without loss of
generality, we may assume the $\a_i$ are distinct.  By Lemma
\ref{zfclem2} choose an $\w$-block $A_i$ of $\a_i$ with $X_i\cap A_i=0$.  By
Lemma \ref{zfclem1} $\bigcap_{i<n}f^{-1}_{\a_i}A_i$ is infinite, and this set
is disjoint from $\bigcup_{i<n}f^{-1}_{\a_i}X_i$. \end{pf}

Let $u$ be any ultrafilter on $\w$ with $u\cap I=0$.  If $u$ were a
$\b$-ultrafilter for some $\b<\w_1$, then choose $\a<\w_1$ indecomposable
with
$\a>\b$.  But then since $f_\a:\w\to\a$ there would be $X\sse\a$, $\tp
X\le\b$,
with $f^{-1}_\a X\in u$.  But by definition $f^{-1}_\a X\in I$ and $I\cap
u=0$.  This completes the proof.

We might comment that it is quite easy to show in ZFC that there is a
non-nowhere dense ultrafilter.  The collection of nowhere dense subsets of
\Q\
forms an ideal, and any ultrafilter on \Q\ disjoint from this ideal is
(isomorphic to) a non-nowhere dense ultrafilter.

The next result complements Theorem \ref{CHthm} in showing that the ordinal
ultrafilters and the topologically-defined ultrafilters do not have much to
do
with one another, although both classes begin with the $P$-points.

\thm{incomp}{Assume Martin's Axiom for $\s$-centered partial orderings. Then
there is an ultrafilter on
$\w$ that is discrete but is not a $J_\a$-ultrafilter for any $\a<\w_1$.}

\begin{pf} We assume the construction of the $f_\a:\w\to\a$ and the ideal $I$
of the previous proof.  For this part of the proof assume $\neg$CH.  We deal
with CH later.

Let $F$ be a filter such that $|F|<\c$ and $F\cap I=0$.  Fix $f:\w\to{}^\w2$. 
It will suffice to find $A\sse\w$ such that $\forall X\in F$ $X\cap A\notin
I$
and $f(A)$ is discrete.

Consider the partial ordering $P$ consisting of all $p=(x_p,g_p,T_p)$ such
that
$x_p\in[\w]^{<\w}$, $g_p:x_p\to\w$, if $m,n\in x_p$ and $m\ne n$ then
$[f(m)\rest g_p(m)]\cap[f(n)\rest g_p(n)]=0$, and $T_p$ is a finite set of
sequences of the form $(X,\a,\b,\g,\d)$ where $X\in F$, $\a<\w_1$ is
indecomposable, $\b<\g\le\a$, $\d$ is indecomposable and if
$Z=\bigcup\set{f^{-1}[f(n)\rest g_p(n)]:n\in x_p}$ then $X\cap(\w-Z)\notin I$
and for any $(X,\a,\b,\g,\d)\in T_p$, $f_\a(X\cap(\w-Z))\cap[\b,\g)$ has
order
type $\ge\d$.  We put $p\le q$ iff $x_p\supseteq x_q$, $g_p\supseteq g_q$ and
$T_p\supseteq T_q$.

If $x_p=x_q$ and $g_p=g_q$ then clearly $p$ and $q$ are compatible, so $P$ is
$\s$-centered.

If $G$ is $P$-generic with repect to some dense sets, we let
$A_G=\bigcup\set{x_p:p\in G}$.  We must choose dense sets so that $A_G$ will
serve as the $A$ we wish.  It is clear that $f(A_G)$ will be discrete.

Let $X\in F$ and let $\a<\w_1$ be indecomposable.  It will suffice to show by
induction on indecomposable $\d\le\a$ that for any $\b,\g$ there is a
countable collection of dense sets in $P$ such that if $G$ is generic for
this
collection and $(X,\a,\b,\g,\d)\in \bigcup\set{T_p:p\in G}$ then
$\tp(f(A_G\cap X)\cap[\b,\g))\ge \d$.

First suppose $\d=\w$.  Fix $\b,\g$.  For $n<\w$ let $D_n=\{\,p\in
P:$ either there is no $q\le p$ with $(X,\a,\b,\g,\d)\in T_q$ or else
$(X,\a,\b,\g,\d)\in T_p$ and $f_\a(x_p\cap[\b,\g))$ has cardinality
$\ge n\,\}$.  We check by induction on $n$ that $D_n$ is dense.  Suppose
$p\in
D_n$ and $(X,\a,\b,\g,\d)\in T_p$.  We seek $q\le p$ with $q\in D_{n+1}$. 
Let
$Z=\bigcup\set{f^{-1}[f(k)\rest g_p(k)]: k\in x_p}$ and let $y\sse
X\cap(\w-Z)$
be such that $|y|=|T_p|+2$ and $f_\a y\sse[\b,\g)$.  For each $k\in y$ find
$r_k$
so large that $[f(k)\rest r_k]$ are all disjoint and are disjoint from $Z$. 
Then there is at most one $k$ such that $(\w-Z)-f^{-1}[f(k)\rest r_k]\in I$
and for each $(X',\a',\b',\g',\d')\in T_p$ there is at most one $k$ such that
$f_\a(X'\cap(\w-Z)-f^{-1}[f(k)\rest r_k])\cap[\b',\g')$ has order type
$<\d'$. 
Hence there is at least one $k\in y$ which does not have any of these
properties, and now if $x_q=x_p\cup\{k\}$ and $g_q=g_p\cup \{(k,r_k)\}$ then
$q\le p$ and
$q\in D_{n+1}$.

Now suppose $\d>\w$.  We operate much as before.  If $\d=\w^{\r+1}$ then put
$\a_n=\w^\r$ for all $n$; otherwise let $\a_n$ be a cofinal sequence of
indecomposable ordinals in $\d$.  Let $D_n=\{\,p\in P:$ either it is
impossible that $(X,\a,\b,\g,\d)\in T_p$ or else $(X,\a,\b,\g,\d)\in
T_p$ and there are $\b\le\b_0<\g_0<\b_1<\g_1<\cdots<\b_{n-1}<\g_{n-1}<\g$
such that $(X,\a,\b_i,\g_i,\a_i)\in T_i$ for all $i<n\,\}$, where here
$\b,\g$ are fixed as before.  Suppose now $p\in D_n$ and $(X,\a,\b,\g,\d)\in
T_p$ as before.  Then $(X,\a,\b_i,\g_i,\a_i)\in T_p$ also, for appropriate
choices of $\b_i,\g_i$.  Let $Z=\bigcup\set{f^{-1}[f(k)\rest g_p(k)]:k\in
x_p}$.  Then $f_\a(X\cap(\w-Z))\cap[\b,\g)$ has order type $\ge \d$ so we can
find $\b_n<\g_n<\g$ such that $\g_{n-1}<\b_n$ and
$f_\a(X\cap(\w-Z))\cap[\b_n,\g_n)$ has order type $\ge\a_n$. But now we may
find
$q\le p$ by setting $x_q=x_p$, $g_q=g_p$, and adding $(X,\a,\b_n,\g_n,\a_n)$
to
$T_p$ to get $T_q$, so $q\in D_{n+1}$.  And now, of course, we need to
supplement the $D_n$ by all the dense sets associated with the
$(X,\a,\b,\g,\d)$
for earlier $\d$.  So by induction, if $G$ is generic for all these dense
sets
and $(X,\a,\b,\g,\d)\in T_p$ for $p\in G$ then $f_\a(A_G)\cap[\b,\g)$ will
have
order type $\ge\d$, as desired.  This completes the proof for MA.  For CH we
do
essentially the same thing, but never need more than countably many dense
sets,
and we must dovetail the construction of the $A_G$'s with the definition of
the
$f_\a$.  Details are straightforward. \end{pf}

We conjecture that it is consistent that every non-principal ultrafilter on
$\w$
is not nowhere dense and (simultaneously) is not a $J_\a$-ultrafilter for any
$\a<\w_1$.  Since every $P$-point is an $\w$-ultrafilter (and is
discrete), this would strongly extend Shelah's result that it is consistent
that
there are no
$P$-points.

\end{document}